\documentclass[MSNbibl,nameyear,dvips]{arxbj}
\usepackage{graphicx}

\aid{0}
\volume{19}
\issue{4}
\pubyear{2013}
\firstpage{1327}
\lastpage{1349}
\doi{10.3150/12-BEJSP06} 

\makeatletter

\newcommand{\cal}{\mathcal}

\newcommand{\citecs}[1]{\citeauthor{#1}, \citeyear{#1}}
\renewcommand{\citep}[1]{(\citeauthor{#1}, \citeyear{#1})}

\newcommand{\real}{\mathbb{R}}
\newcommand{\sphere}{\mathbb{S}}

\newcommand{\ang}{\theta}
\newcommand{\dd}{\mathrm{d}}

\newcommand{\bnd}{b_{n,d}}
\newcommand{\bnone}{b_{n,1}}
\newcommand{\bntwo}{b_{n,2}}

\newtheorem{theorem}{Theorem}
\newtheorem{corollary}{Corollary}
\newtheorem{lemma}{Lemma}

\newremark{example}{Example}

\makeatother

\begin{document}
\begin{frontmatter}

\title{Strictly and non-strictly positive definite functions on spheres}
\runtitle{Positive definite functions on spheres}

\begin{aug}
\author{\fnms{Tilmann} \snm{Gneiting}\corref{}\ead[label=e1]{t.gneiting@uni-heidelberg.de}\ead[label=e2,url]{www.math.uni-heidelberg.de/spatial}}
\runauthor{T. Gneiting} 
\address{Institut f\"ur Angewandte Mathematik, Universit\"at
Heidelberg, Im Neuenheimer Feld 294, D-69120 Heidelberg,
Germany. \printead{e1}}
\end{aug}


%
\begin{abstract}
Isotropic positive definite functions on spheres play important roles
in spatial statistics, where they occur as the correlation functions
of homogeneous random fields and star-shaped random particles. In
approximation theory, strictly positive definite functions serve as
radial basis functions for interpolating scattered data on spherical
domains. We review characterizations of positive definite functions
on spheres in terms of Gegenbauer expansions and apply them to
dimension walks, where monotonicity properties of the Gegenbauer
coefficients guarantee positive definiteness in higher dimensions.
Subject to a natural support condition, isotropic positive definite
functions on the Euclidean space $\real^3$, such as Askey's and
Wendland's functions, allow for the direct substitution of the
Euclidean distance by the great circle distance on a one-, two- or
three-dimensional sphere, as opposed to the traditional approach,
where the distances are transformed into each other. Completely
monotone functions are positive definite on spheres of any dimension
and provide rich parametric classes of such functions, including
members of the powered exponential, Mat\'ern, generalized Cauchy and
Dagum families. The sine power family permits a continuous
parameterization of the roughness of the sample paths of a Gaussian
process. A collection of research problems provides challenges for
future work in mathematical analysis, probability theory and spatial
statistics.
\end{abstract}

%

\pdfkeywords{completely monotone,
covariance localization, fractal index,
interpolation of scattered data,
isotropic, locally supported,
multiquadric, Polya criterion,
radial basis function,
Schoenberg coefficients}

\end{frontmatter}


\section{Introduction} \label{secintroduction}

Recently, there has been renewed interest in the study of positive
definite functions on spheres, motivated in part by applications in
spatial statistics, where these functions play crucial roles as the
correlation functions of random fields on spheres
(\citecs{Banerjee2005}; \citecs{Huang2011}), including the case of star-shaped
random particles
\citep{Hansen2011}. In a related development in approximation theory,
strictly positive definite functions arise as radial basis functions
for interpolating scattered data on spherical domains
(\citecs{Fasshauer1998}; \citecs{Cavoretto2010}; \citecs{LeGia2010}).

Specifically, let $k \geq2$ be an integer, and let $\sphere^{k-1} =
\{ x \in\real^k: \| x \| = 1 \}$ denote the unit sphere in the
Euclidean space $\real^k$, where we write $\| x \|$ for the Euclidean
norm of $x \in\real^k$. A function $h: \sphere^d \times\sphere^d
\to\real$ is \textit{positive definite} if
%
%
\begin{equation}
\label{eqpd} \sum_{i=1}^n \sum
_{j=1}^n c_i c_j
h(x_i, x_j) \geq0
\end{equation}
for all finite systems of pairwise distinct points $x_1,\ldots, x_n
\in\sphere^d$ and constants $c_1,\ldots, c_n \in\real$. A positive
definite function $h$ is \textit{strictly} positive definite if the
inequality in (\ref{eqpd}) is strict, unless $c_1 = \cdots= c_n =
0$, and is \textit{non-strictly} positive definite otherwise. The
function $h: \sphere^d \times\sphere^d \to\real$ is \textit{isotropic}
if there exists a function $\psi: [0,\pi] \to\real$ such that
%
%
\begin{equation}
\label{eqh} h(x,y) = \psi\bigl(\ang(x,y)\bigr) \quad\mbox{for all}
\quad x, y
\in\sphere^d,
\end{equation}
where $\ang(x,y) =
\arccos( \langle x, y \rangle)$ is the great circle, spherical or
geodesic distance on $\sphere^d$, and $\langle\cdot, \cdot\rangle$
denotes the scalar or inner product in $\real^{d+1}$. Thus, an
isotropic function on a sphere depends on its arguments via the great circle
distance $\theta(x,y)$ or, equivalently, via the inner product
$\langle
x, y \rangle$ only.

For $d = 1, 2,\ldots,$ we write $\Psi_d$ for the class of the
continuous functions $\psi: [0,\pi] \to\real$ with $\psi(0) = 1$
which are such that the associated isotropic function $h$ in
(\ref{eqh}) is positive definite. Hence, we may identify the class
$\Psi_d$ with the correlation functions of the mean-square continuous,
stationary and isotropic random fields on the sphere $\sphere^d$
\citep{Jones1963}. We write $\Psi_d^+$ for the class of the
continuous functions $\psi: [0,\pi] \to\real$ with $\psi(0) = 1$
such that the isotropic function $h$ in (\ref{eqh}) is strictly
positive definite.

We study the classes $\Psi_d$ and $\Psi_d^+$ with particular attention
to the practically most relevant case of a two-dimensional sphere,
such as planet Earth. Applications abound, with atmospheric data
assimilation and the reconstruction of the global temperature and
green house gas concentration records being key examples. Section
\ref{secEuclidean} summarizes related results for isotropic positive
definite functions on Euclidean spaces and studies their restrictions
to spheres. In Section~\ref{secspheres}, we review characterizations
of strictly and non-strictly positive definite functions on spheres in
terms of expansions in ultraspherical or Gegenbauer polynomials.
Monotonicity properties of the Gegenbauer coefficients in any given
dimension guarantee positive definiteness in higher dimensions.

%
%
\begin{table}
\caption{Parametric families of correlation functions and radial basis
functions on one-, two- and three-dimensional spheres in terms of
the great circle distance, $\theta\in[0,\pi]$. Here $c$ is a
scale or support parameter, $\delta$ and $\tau$ are shape
parameters, $\alpha$ or $\nu$ is a smoothness parameter, and we
write $t_+ = \max(t,0)$ for $t \in\real$. The parameter range
indicated guarantees membership in the classes $\Psi_d^+$, where $d
= 1, 2$ and $3$. For details, see Section
\protect\ref{secapplications}}
\label{tabPsi2}
\begin{tabular*}{\tablewidth}{@{\extracolsep{\fill}}lll@{}}
\hline
Family
& Analytic expression
& Parameter range \\
\hline
Powered exponential
& $\psi(\theta) = \exp( - ( \frac{\theta}{c}
)^{
\alpha} )$
& $c > 0; \alpha\in(0,1]$ \\[3pt]
Mat\'ern & $\psi(\theta) = \frac{2^{\nu-1}}{\Gamma
(\nu)}
( \frac{\theta}{c} )^{ \nu} K_\nu( \frac
{\theta
}{c} )$
& $c > 0; \nu\in(0,\frac{1}{2}]$ \\[3pt]
Generalized Cauchy & $\psi(\theta) = ( 1 +
( \frac{\theta}{c} )^{ \alpha})^{- \tau
/\alpha}$
& $c > 0; \tau> 0; \alpha\in(0,1]$ \\[3pt]
Dagum & $\psi(\theta) = 1 - (
( \frac{\theta}{c} )^{ \tau} / ( 1 +
\frac
{\theta}{c} )^{ \tau}
)^{\alpha/\tau}$
& $c > 0; \tau\in(0,1]; \alpha\in(0,\tau)$ \\[3pt]
Multiquadric
& $\psi(\theta) = (1-\delta)^{2\tau} / (1 + \delta^2 - 2
\delta\cos\theta)^\tau$
& $\tau> 0; \delta\in(0,1)$ \\[3pt]
Sine power
& $\psi(\theta) = 1 - ( \sin\frac{\theta}{2} )^\alpha$
& $\alpha\in(0,2)$ \\[3pt]
Spherical & $\psi(\theta) = ( 1 + \frac{1}{2}
\frac
{\theta}{c} )
( 1 - \frac{\theta}{c} )_+^2$ & $c > 0$ \\[3pt]
Askey & $\psi(\theta) = ( 1 - \frac{\theta}{c}
)^\tau_+$
& $c > 0; \tau\geq2$ \\[3pt]
$C^2$-Wendland & $\psi(\theta) = ( 1 + \tau
\frac{\theta}{c} )
( 1 - \frac{\theta}{c} )^\tau_+$ & $c \in(0,\pi];
\tau
\geq4$ \\[3pt]
$C^4$-Wendland & $\psi(t) = ( 1 + \tau
\frac
{\theta}{c} +
\frac{\tau^2-1}{3} \frac{\theta^2}{c^2} ) ( 1 -
\frac
{\theta}{c} )^\tau_+$
& $c \in(0,\pi]; \tau\geq6$ \\
\hline
\end{tabular*}
\end{table}
%

The core of the paper is Section~\ref{secapplications}, where we
supply easily applicable conditions for membership in the classes
$\Psi_d$ and $\Psi_d^+$, and use them to derive rich parametric
families of correlation functions and radial basis functions on
spheres. In particular, compactly supported, isotropic positive
definite functions on the Euclidean space $\real^3$ remain positive
definite with the great circle distance on a one-, two- or
three-dimensional sphere as argument. Consequently, compactly
supported radial basis functions on $\real^3$, such as Askey's and
Wendland's functions, translate directly into locally supported radial
basis functions on the sphere. Criteria of P\'olya type guarantee
membership in the classes $\Psi_d^+$, and we supplement recent results
of \citet{Beatson2011} in the case of spheres of dimension $d \leq7$.
Completely monotone functions are strictly positive definite on
spheres of any dimension, and include members of the powered
exponential, Mat\'ern, generalized Cauchy and Dagum families. Some of
the key results in the cases of one-, two- and three-dimensional
spheres are summarized in Table~\ref{tabPsi2}.

The paper ends in Section~\ref{secdiscussion}, where we refer to a
collection of research problems in mathematical analysis, probability
theory and statistics that offer challenges of diverse difficulty and
scope, with details given in a supplemental article
\citep{Gneiting2012}.


\section{Isotropic positive definite functions on Euclidean spaces}
\label{secEuclidean}

In this expository section, we review basic results about isotropic
positive definite functions on Euclidean spaces and study their
restrictions to spheres.

Recall that a function $h: \real^d \times\real^d \to\real$ is
\textit{positive definite} if the inequality (\ref{eqpd}) holds for all
finite systems of pairwise distinct points $x_1,\ldots, x_n \in
\real^d$ and constants $c_1,\ldots, c_n \in\real$. It is \textit
{strictly positive definite} if the inequality is strict unless $c_1
= \cdots= c_n = 0$. The function $h: \real^d \times\real^d \to
\real$ is \textit{radial}, \textit{spherically symmetric} or \textit
{isotropic} if there exists a function $\varphi: [0,\infty) \to\real$
such that
%
%
\begin{equation}
\label{eqhrad} h(x,y) = \varphi\bigl(\| x - y \|\bigr) \quad\mbox{for all} \quad
x, y
\in\real^d.
\end{equation}
For an integer $d \geq1,$ we denote by $\Phi_d$ the class of the
continuous functions $\varphi: [0,\infty) \to\real$ with $\varphi(0)
= 1$ such that the function $h$ in (\ref{eqhrad}) is positive
definite. Thus, we may identify the members of the class $\Phi_d$
with the characteristic functions of spherically symmetric probability
distributions, or with the correlation functions of mean-square
continuous, stationary and isotropic random fields on $\real^d$.

\citet{Schoenberg1938a} showed that a function $\varphi: [0,\infty)
\to\real$ belongs to the class $\Phi_d$ if, and only if, it is of the
form
%
%
\begin{equation}
\label{eqPhid} \varphi(t) = \int_{[0,\infty)}
\Omega_d(rt) \,\dd F(r) \quad\mbox{for} \quad t \geq0,
\end{equation}
where $F$ is a uniquely determined probability measure on $[0,\infty)$
and
\[
\Omega_d(t) = \Gamma(d/2) \biggl( \frac{2}{t}
\biggr)^{ (d-2)/2} J_{(d-2)/2}(t),
\]
with $J$ being a Bessel function (\citecs{DLMF}, Section 10.2). The
classes $\Phi_d$ are nonincreasing in $d$,
\[
\Phi_1 \supset\Phi_2 \supset\cdots\supset
\Phi_\infty= \bigcap_{d=1}^\infty
\Phi_d,
\]
with the inclusions being strict. If $d \geq2$, any nonconstant
member $\varphi$ of the class $\Phi_d$ corresponds to a strictly
positive definite function on $\real^d$ (\citecs{Sun1993},
Theorem 3.8)
and thus it belongs to the class $\Phi_d^+$. Finally, as shown by
\citet{Schoenberg1938a}, the class $\Phi_\infty$ consists of the
functions $\varphi$ of the form
%
%
\begin{equation}
\label{eqSchoenberg} \varphi(t) = \int_{[0,\infty)} \exp
\bigl(-r^2 t^2\bigr) \,\dd F(r) \quad\mbox{for} \quad t
\geq0,
\end{equation}
where $F$ is a uniquely determined probability measure on $[0,\infty)$.

From the point of view of isotropic random fields, the asymptotic
decay of the function $\varphi$ at the origin determines the
smoothness of the associated Gaussian sample path
\citep{Adler1981}. Specifically, if
%
%
\begin{equation}
\label{eqfractalphi} \varphi(0) - \varphi(t) = {\cal O}\bigl(t^\alpha
\bigr) \quad\mbox{as} \quad t \downarrow0
\end{equation}
for some $\alpha\in(0,2]$, which we refer to as the \textit{fractal
index}, the graph of a Gaussian sample path has fractal or Hausdorff
dimension $D = d + 1 - \frac{\alpha}{2}$ almost surely. If $\alpha=
2$, the sample
path is smooth and differentiable and its Hausdorff dimension, $D$,
equals its topological dimension,~$d$. If $\alpha\in(0,2)$, the
sample path is non-differentiable. Table~\ref{tabPhiinfty} shows
some well known parametric families within the class $\Phi_\infty$,
namely the powered exponential family \citep{Yaglom1987}, the Mat\'ern
class \citep{GuttorpGneiting2006}, the generalized Cauchy family
\citep{GneitingSchlather2004}, and the Dagum class \citep{Berg2008},
along with the corresponding fractal indices. Table~\ref{tabcompact}
shows compactly supported families within the class $\Phi_3$ that have
been introduced and studied by \citet{Wendland1995} and
Gneiting (\citeyear{Gneiting1999a,Gneiting2002}).

%
%
\begin{table}
\caption{Parametric families of members of the class $\Phi_\infty$.
Here $K_\nu$ denotes the modified Bessel function of the second
kind of order $\nu$ (\citecs{DLMF}, Section 10.2), $c$ is a scale
parameter, $\alpha$ and $\nu$ are smoothness parameters, and $\tau$
is a shape parameter, respectively}
\label{tabPhiinfty}
\begin{tabular*}{\tablewidth}{@{\extracolsep{\fill}}llll@{}}
\hline
Family
& Analytic expression
& Parameter range
& Fractal index \\
\hline
$
\begin{array}{@{}l} \mbox{Powered} \\ \mbox{exponential}
\end{array}
$&
$\varphi(t) = \exp( - ( \frac{t}{c} )^{
\alpha})$
& $c > 0; \alpha\in(0,2]$ & $\alpha\in(0,2]$ \\[4pt]
$
\begin{array}{@{}l} \mbox{Mat\'ern}
\end{array}
$ &
$\varphi(t) = \frac{2^{\nu-1}}{\Gamma(\nu)}
( \frac{t}{c} )^{ \nu} K_\nu( \frac{t}{c}
)$
& $c > 0; \nu> 0$ & $\min(2\nu,2) \in(0,2]$ \\[4pt]
$
\begin{array}{@{}l} \mbox{Generalized} \\ \mbox{Cauchy}
\end{array}
$ & $\varphi(t) = ( 1 + ( \frac{t}{c} )^{\alpha})^{-
\tau/\alpha}$
& $c > 0; \alpha\in(0,2]; \tau> 0$ & $\alpha\in(0,2]$ \\[4pt]
$
\begin{array}{@{}l} \mbox{Dagum}
\end{array}
$ &
$\varphi(t) = 1 - ( ( \frac{t}{c} )^{ \tau}
/ ( 1 + \frac{t}{c} )^{ \tau}
)^{ \alpha/\tau}$
& $c > 0; \tau\in(0,2]; \alpha\in(0,\tau)$ & $\alpha\in
(0,\tau)$ \\
\hline
\end{tabular*}
\end{table}
%

%
%
\begin{table}[b]
\tablewidth=335pt
\caption{Parametric families of compactly supported members of the
class $\Phi_3$ with support parameter $c$ and shape parameter
$\tau$}
\label{tabcompact}
\begin{tabular*}{\tablewidth}{@{\extracolsep{\fill}}llll@{}}
\hline
Family
& Analytic expression
& Parameter range
& Fractal index \\
\hline
Spherical
& $\varphi(t) = ( 1 + \frac{1}{2} \frac{t}{c} )
( 1
- \frac{t}{c} )_+^2$
& $c > 0$ & 1 \\[3pt]
Askey
& $\varphi(t) = ( 1 - \frac{t}{c} )^\tau_+$ & $c > 0;
\tau\geq2$ & 1 \\[3pt]
$C^2$-Wendland
& $\varphi(t) = ( 1 + \tau\frac{t}{c} ) (
1 -
\frac{t}{c} )^\tau_+$
& $c > 0; \tau\geq4$ & 2 \\[3pt]
$C^4$-Wendland
& $\varphi(t) = ( 1 + \tau\frac{t}{c} + \frac{\tau^2-1}{3}
\frac{t^2}{c^2} )
( 1 - \frac{t}{c} )^\tau_+$
& $c > 0; \tau\geq6$ & 2 \\
\hline
\end{tabular*}
\end{table}
%

\citet{Yadrenko1983} pointed out that if $\varphi$ is a member of the
class $\Phi_k$ for some integer $k \geq2$, then the function defined
by
%
%
\begin{equation}
\label{eqYadrenko} \psi(\theta) = \varphi\biggl( 2 \sin\frac{\theta}{2}
\biggr) \quad\mbox{for} \quad\theta\in[0,\pi]
\end{equation}
corresponds to the restriction of the isotropic function $h(x_i, x_j)
= \varphi(\| x_i - x_j \|)$ from $\real^k \times\real^k$ to
$\sphere^{k-1} \times\sphere^{k-1}$, with the chordal or Euclidean
distance expressed in terms of the great circle dis\-tance, $\theta$,
on the sphere $\sphere^{k-1}$. Thus, given any nonconstant member
$\varphi$ of the class $\Phi_k$, the function $\psi$ defined by
(\ref{eqYadrenko}) belongs to the class $\Psi_{k-1}^+$.

Various authors have argued in favor of this construction, including
Fasshauer and Schumaker (\citeyear{Fasshauer1998}), \citet{Gneiting1999a},
\citet{NarcowichWard2002} and \citet{Banerjee2005}, as it readily
generates parametric families of isotropic, strictly positive definite
functions on spheres and retains the interpretation of scale, support,
shape and smoothness parameters. In particular, the mapping
(\ref{eqYadrenko}) from $\varphi\in\Phi_k$ to $\psi\in\Psi_{k-1}$
preserves the fractal index, in the sense that
%
%
\begin{equation}
\label{eqfractalpsi} \psi(0) - \psi(\theta) = {\cal O}\bigl(
\theta^\alpha\bigr) \quad\mbox{as} \quad\theta\downarrow0
\end{equation}
if $\varphi$ has fractal index $\alpha\in(0,2]$, as defined in
(\ref{eqfractalphi}). Nevertheless, the approach is of limited
flexibility. For example, if $\varphi\in\Phi_3$ it is readily seen
from (\ref{eqPhid}) that the function $\psi$ defined in
(\ref{eqYadrenko}) does not admit values less than $\inf_{t > 0}
t^{-1} \sin t = - 0.2127\ldots\,$. Furthermore, the mapping
$\theta\mapsto2 \sin\frac{\theta}{2}$ in the argument of
$\varphi$ in (\ref{eqYadrenko}), while being essentially linear for
small $\theta$, is counter to spherical geometry for larger values of
the great circle distance, and thus may result in physically
unrealistic distortions. In this light, we now turn to characterizations
and constructions of positive definite functions that operate directly
on a sphere.


\section{Isotropic positive definite functions on spheres} \label{secspheres}

Let $d \geq1$ be an integer. As defined in Section
\ref{secintroduction}, the classes $\Psi_d$ and $\Psi_d^+$ consist of
the functions $\psi: [0,\pi] \to\real$ with $\psi(0) = 1$ which are
such that the isotropic function $h$ in (\ref{eqh}) is positive
definite or strictly positive definite, respectively. Furthermore, we
consider the class $\Psi_d^- = \Psi_d \setminus\Psi_d^+$ of the
non-strictly positive definite functions, and we define
\[
\Psi_\infty= \bigcap_{d =1}^\infty
\Psi_d, \qquad\Psi_\infty^+ = \bigcap
_{d =1}^\infty\Psi_d^+ \qquad\mbox{and}
\qquad\Psi_\infty^- = \bigcap_{d =1}^\infty
\Psi_d^-.
\]
\citet{Schoenberg1942} noted that the classes $\Psi_d$ and
$\Psi_\infty$ are convex, closed under products, and closed under
limits, provided the limit function is continuous. The classes
$\Psi_d^+$ and $\Psi_\infty^+$ are convex and closed under products,
but not under limits.

We proceed to review characterizations in terms of ultra\-spherical or
Gegenbauer expansions, for which we require classical results on
orthogonal polynomials (Digital Library of Mathematical
Functions (\citeyear{DLMF}), Section 18.3). Given $\lambda>
0$ and an integer $n \geq0$, the function $C_n^\lambda(\cos\theta)$
is defined by the expansion
%
%
\begin{equation}
\label{eqGegenbauer} \frac{1}{ ( 1 + r^2 - 2r \cos\theta)^\lambda} =
\sum_{n=0}^\infty
r^n C_n^{\lambda}(\cos\theta) \quad\mbox{for}
\quad\theta\in[0,\pi],
\end{equation}
where $r \in(-1,1)$, and $C_n^\lambda$ is the ultraspherical or
Gegenbauer polynomial of degree $n$, which is even if $n$ is even and
odd if $n$ is odd. For reference later on, we note that
$C_n^{\lambda}(1) = \Gamma(n + 2\lambda)/(n!
\Gamma(2\lambda))$. If $\lambda= 0$, we follow
\citet{Schoenberg1942} and set
\[
C_n^0(\cos\theta) = \cos(n\theta) \quad\mbox{for} \quad
\theta\in[0,\pi].
\]
By equation (3.42) of \citet{Askey1969},
%
%
\begin{equation}
\label{eqGegenbauerrel} C_n^\lambda(\cos\theta) =
\frac{\Gamma(\nu)}{\Gamma(\lambda) \Gamma(\lambda-\nu)} \sum_{k
= 0}^{[n/2]}
\frac{(n-2k+\nu) \Gamma(k+\lambda-\nu) \Gamma
(n-k+\lambda
)}{k! \Gamma(n-k+\nu+1)} C_{n-2k}^\nu(\cos\theta)
\end{equation}
whenever $\lambda> \nu\geq0$ and $\theta\in[0,\pi]$, where the
coefficients on the right-hand side are strictly positive. This
classical result of Gegenbauer will be used repeatedly in the sequel.

The following theorem summarizes Schoenberg's (\citeyear{Schoenberg1942}) classical
characterization of the classes $\Psi_d$ and $\Psi_\infty$, along with
more recent results of \citet{Menegatto1994}, Chen, Menegatto and Sun (\citeyear{Chen2003}) and
\citet{Menegatto2006} for the classes $\Psi_d^+$ and $\Psi_\infty^+$.

%
\begin{theorem} \label{thPsi}
Let $d \geq1$ be an integer.
\begin{itemize}[(a)]
\item[(a)]
The class $\Psi_d$ consists of the functions of the form
%
%
\begin{equation}
\label{eqPsid} \psi(\theta) = \sum_{n = 0}^\infty
\bnd\frac{C_n^{(d-1)/2}(\cos
\theta)}{C_n^{(d-1)/2}(1)} \quad\mbox{for} \quad\theta\in[0,\pi],
\end{equation}
with nonnegative, uniquely determined coefficients $\bnd$ such
that $\sum_{n = 0}^\infty\bnd= 1$. If $d \geq2$, the class
$\Psi_d^+$ consists of the functions in $\Psi_d$ with the
coefficients $\bnd$ being strictly positive for infinitely many even
and infinitely many odd integers $n$. The class $\Psi_1^+$
consists of the functions in $\Psi_1$ which are such that, given any
two integers $0 \leq j < n$, there exists an integer $k \geq0$
with the coefficient $b_{j+kn,1}$ being strictly positive.
\item[(b)]
The class $\Psi_\infty$ consists of the functions of the form
%
%
\begin{equation}
\label{eqPsiinfty} \psi(\theta) = \sum
_{n = 0}^\infty b_n (\cos\theta)^n
\quad\mbox{for} \quad\theta\in[0,\pi],
\end{equation}
with nonnegative, uniquely determined coefficients $b_n$ such that
$\sum_{n = 0}^\infty b_n = 1$. The class $\Psi_\infty^+$ consists
of the functions in $\Psi_\infty$ with the coefficients $b_n$
being strictly positive for infinitely many even and infinitely many
odd integers $n$.
\end{itemize}
\end{theorem}

If $d = 1$, Schoenberg's representation (\ref{eqPsid}) reduces to
the general form,
%
%
\begin{equation}
\label{eqPsi1} \psi(\theta) = \sum_{n = 0}^\infty
\bnone\cos(n\theta) \quad\mbox{for} \quad\theta\in[0,\pi],
\end{equation}
of a member of the class $\Psi_1$, that is, a standardized, continuous
positive definite function on the circle. By basic Fourier calculus,
%
%
\begin{equation}
\label{eqb1} b_{0,1} = \frac{1}{\pi} \int_0^\pi
\psi(\theta) \,\dd\theta\qquad\mbox{and} \qquad\bnone= \frac{2}{\pi}
\int
_0^\pi\cos(n\theta) \psi(\theta) \,\dd\theta
\end{equation}
for integers $n \geq1$. If $d = 2$, the representation
(\ref{eqPsid}) yields the general form,
%
%
\begin{equation}
\label{eqPsi2} \psi(\theta) = \sum_{n = 0}^\infty
\frac{\bntwo}{n+1} P_n(\cos\theta) \quad\mbox{for} \quad\theta
\in[0,\pi],
\end{equation}
of a member of the class $\Psi_2$ in terms of the Legendre polynomial
$P_n$ of integer order $n \geq0$. As \citet{Bingham1973} noted, the
classes considered here are convex, and so (\ref{eqPsid}),
(\ref{eqPsiinfty}), (\ref{eqPsi1}) and (\ref{eqPsi2}) can be
interpreted as Choquet representations in terms of extremal members,
which are non-strictly positive definite.

The general form (\ref{eqPsiinfty}) of the members of the class
$\Psi_\infty$ can be viewed as a power series in the variable $\cos
\theta$ with nonnegative coefficients. To give an example, a standard
Taylor expansion shows that the function defined by
%
%
\begin{equation}
\label{eqmultiquadric} \psi(\theta) = \frac{(1-\delta)^{2\tau}}{(1 +
\delta^2 - 2 \delta
\cos
\theta)^\tau} \quad\mbox{for} \quad
\theta\in[0,\pi]
\end{equation}
admits such a representation if $\tau> 0$ and $\delta\in(0,1)$;
furthermore, the coefficients are strictly positive, whence $\psi$
belongs to the class $\Psi_\infty^+$. When $d = 2$, the special cases
in (\ref{eqmultiquadric}) with the parameter $\tau$ fixed at
$\frac{1}{2}$ and $\frac{3}{2}$ have been known as the inverse
multiquadric and the Poisson spline, respectively
\citep{Cavoretto2010}. In this light, we refer to
(\ref{eqmultiquadric}) as the \textit{multiquadric} family.

The \textit{sine power} function of \citet{Soubeyrand2008},
%
%
\begin{equation}
\label{eqsinepower} \psi(\theta) = 1 - \biggl( \sin\frac{\theta}{2}
\biggr)^{
\alpha} \quad\mbox{for} \quad\theta\in[0,\pi],
\end{equation}
also admits the representation (\ref{eqPsiinfty}) if $\alpha\in
(0,2]$. Moreover, the coefficients are strictly positive if $\alpha
\in(0,2)$, whence $\psi$ belongs to the class $\Psi_\infty^+$. The
parameter $\alpha$ corresponds to the fractal index in the
relationship (\ref{eqfractalpsi}) and thus parameterizes the
roughness of the sample paths of the associated Gaussian
process.

The classes $\Psi_d$ are nonincreasing in $d \geq1$, with the
following result revealing details about their structure. The
statement in part (a) might be surprising, in that a function that is
non-strictly positive definite on the sphere $\sphere^d$ cannot be
strictly positive definite on a lower-dimensional sphere, including
the case $\sphere^1$ of the circle.

%
\begin{corollary} \label{corPsi}
\mbox{}
\begin{itemize}[(a)]
\item[(a)]
If $\psi\in\Psi_d^+ \cap\Psi_{d'}$ for positive integers $d$
and $d'$, then $\psi\in\Psi_{d'}^+$. Similarly, if $\psi\in
\Psi_d^- \cap\Psi_{d'}$ for positive integers $d$ and $d'$,
then $\psi\in\Psi_{d'}^-$. In particular,
%
%
\begin{equation}
\label{eqPsiinfty+-} \Psi_\infty= \Psi_\infty^+ \cup
\Psi_\infty^-,
\end{equation}
where the union is disjoint.
\item[(b)]
The classes $\Psi_d^+$ are nonincreasing in $d$,
\[
\Psi_1^+ \supset\Psi_2^+ \supset\cdots\supset
\Psi_\infty^+ = \bigcap_{d \geq1}
\Psi_d^+,
\]
with the inclusions being strict.
\item[(c)]
The classes $\Psi_d^-$ are nonincreasing in $d$,
\[
\Psi_1^- \supset\Psi_2^- \supset\cdots\supset
\Psi_\infty^- = \bigcap_{d \geq1}
\Psi_d^-,
\]
with the inclusions being strict.\vadjust{\goodbreak}
\end{itemize}
\end{corollary}

\begin{pf}
In part (a) we use an argument of \citet{Narcowich1995},
and our key tool in proving that the inclusions in parts (b) and (c)
are strict is Gegenbauer's relationship (\ref{eqGegenbauerrel}).
\begin{itemize}[(a)]
\item[(a)]
If $d \geq d'$, it is trivially true that $\psi\in\Psi_d^+$ implies
$\psi\in\Psi_{d'}^+$. If $d < d'$, suppose, for a contradiction,
that $\psi\in\Psi_d^+ \cap\Psi_{d'}^-$. By part (a) of Theorem
\ref{thPsi} applied in dimension $d' \geq2$, $\psi$ is either an
even function plus an odd polynomial in $\cos\theta$, or an odd
function plus an even polynomial in $\cos\theta$. By Theorem 2.2 of
\citet{Menegatto1995}, we conclude that $\psi\in\Psi_1^-$, for the
desired contradiction to the assumption that $\psi\in\Psi_d^+$. The
proof of the second claim is analogous, and the statement in
(\ref{eqPsiinfty+-}) then is immediate.
\item[(b)]
The inclusion $\Psi_{d+1}^+ \subseteq\Psi_d^+$ is trivially true. To
see that the inclusion is strict, note that by part (a) of Theorem
\ref{thPsi} and the relationship (\ref{eqGegenbauerrel}) any $\psi
\in\Psi_{d+1}^+$ admits the representation (\ref{eqPsid}) with
$\bnd$ being strictly positive for all $n \geq0$. However, there are
members of the class $\Psi_d^+$ with $\bnd= 0$ for at least one
integer $n \geq0$, which thus do not belong to the class $\Psi_{d+1}^+$.
\item[(c)]
The inclusion $\Psi_{d+1}^- \subseteq\Psi_d^-$ is immediate from part
(a). To demonstrate that the inclusion is strict, we show that if $n
\geq2$ then the function $\psi(\theta) = C_n^{(d-1)/2}(\cos\theta)$
belongs to $\Psi_d^-$ but not to $\Psi_{d+1}^-$. For a contradiction,
suppose that $\psi\in\Psi_{d+1}$. Then by the relationship
(\ref{eqGegenbauerrel}), $\psi$ is of the form (\ref{eqPsid}) with
at least two distinct coefficients $\bnd$ being\vspace*{2pt} strictly positive,
contrary to the uniqueness of the Gegenbauer coefficients.\qed
\end{itemize}
\noqed\end{pf}

The following classical result is a consequence of Schoenberg's (\citeyear{Schoenberg1942})
representation (\ref{eqPsid}) and the orthogonality properties of
the ultraspherical or Gegenbauer polynomials.

%
\begin{corollary} \label{corSchoenberg}
Let $d \geq2$ be an integer. A continuous function $\psi:
[0,\pi] \to\real$ with $\psi(0) = 1$ belongs to the class
$\Psi_d$ if and only if
%
%
\begin{equation}
\label{eqbnd} \bnd= \frac{2n+d-1}{2^{3-d} \pi} \frac{(\Gamma(\frac
{d-1}{2}))^2}{\Gamma(d-1)} \int
_0^\pi C_n^{(d-1)/2}(\cos\theta)
(\sin\theta)^{d-1} \psi(\theta) \,\dd\theta\geq0
\end{equation}
for all integers $n \geq0$. Furthermore, the coefficient $\bnd$
in the Gegenbauer expansion (\ref{eqPsid}) of a function $\psi\in
\Psi_d$ equals the above value.
\end{corollary}

Interesting related results include Theorem 4.1 of
\citet{NarcowichWard2002}, which expresses the Gegenbauer coefficient
$\bnd$ of a function $\psi$ of Yadrenko's form (\ref{eqYadrenko}),
where $k = d + 1$, in terms of the Fourier transform of the respective
member $\varphi$ of the class $\Phi_{d+1}$, Theorem 6.2 of
\citet{LeGia2010}, which deduces asymptotic estimates, and the recent
findings of \citet{Ziegel2012} on convolution roots and smoothness
properties of isotropic positive definite functions on spheres.

Adopting terminology introduced by \citet{Daley2012} and
\citet{Ziegel2012}, we refer to the Gegenbauer coefficient $\bnd$ as a
\textit{$d$-Schoenberg coefficient} and to the sequence $(b_{n,d})_{n =
0, 1, 2, \ldots}$ as a \textit{$d$-Schoenberg sequence}. We now provide
formulas that express the $d$-Schoenberg coefficient $b_{n,d+2}$ in
terms of $\bnd$ and $b_{n+2,d}$. In special cases, closely related
results were used by \citet{Beatson2011}.

%
\begin{corollary} \label{corbnd}
Consider the $d$-Schoenberg coefficients $\bnd$ in the Gegenbauer
expansion (\ref{eqPsid}) of the members of the class $\Psi_d$.
\begin{itemize}[(a)]
\item[(a)]
It is true that
%
%
\begin{equation}
\label{eqdw130} b_{0,3} = b_{0,1} - \tfrac{1}{2}
b_{2,1}
\end{equation}
and
%
%
\begin{equation}
\label{eqdw13} b_{n,3} = \tfrac{1}{2} (n+1) ( b_{n,1} -
b_{n+2,1} )
\end{equation}
for all integers $n \geq1$.
\item[(b)]
If $d \geq2$, then
%
%
\begin{equation}
\label{eqdw} b_{n,d+2} = \frac{(n+d-1) (n+d)}{d (2n+d-1)} \bnd- \frac
{(n+1) (n+2)}{d (2n+d+3)}
b_{n+2,d}
\end{equation}
for all integers $n \geq0$.
\end{itemize}
\end{corollary}

\begin{pf}
We take advantage of well known recurrence relations for
trigonometric functions and Gegenbauer polynomials.
\begin{itemize}[(a)]
\item[(a)]
The identity in (\ref{eqdw130}) holds true because $2 (\sin
\theta
)^2 =
1 - \cos(2\theta)$ for $\theta\in[0,\pi]$, so that
\[
b_{0,3} = \frac{2}{\pi} \int_0^\pi(
\sin\theta)^2 \psi(\theta) \,\dd\theta= \frac{1}{\pi} \int
_0^\pi\bigl(1 - \cos(2\theta)\bigr) \psi(\theta)
\,\dd\theta= b_{0,1} - \frac{1}{2} b_{2,1}.
\]
To establish (\ref{eqdw13}), we note that
\[
C_n^1(\cos\theta) (\sin\theta)^2 = \sin
\theta\sin\bigl((n+1)\theta\bigr) = \tfrac{1}{2} \bigl( \cos(n\theta) -
\cos
\bigl((n+2)\theta\bigr) \bigr)
\]
for $\theta\in[0,\pi]$. In view of the formulas (\ref{eqbnd}) and
(\ref{eqb1}) for the Gegenbauer coefficients $b_{n,3}$ and the
Fourier cosine coefficients $\bnone$, respectively, the claim follows
easily if we integrate the above equality over $\theta\in[0,\pi]$.
\item[(b)]
By equation (18.9.8) of \citet{DLMF},
\begin{eqnarray*}
C_n^{(d+1)/2}(\cos\theta) (\sin\theta)^2 & = &
\frac{(n+d-1) (n+d)}{(d-1) (2n+d+1)} C_n^{(d-1)/2}(\cos\theta)
\\
&&{} - \frac{(n+1) (n+2)}{(d-1) (2n+d+1)} C_{n+2}^{(d-1)/2}(\cos
\theta).
\end{eqnarray*}
In view of the formula (\ref{eqbnd}) for the $d$-Schoenberg
coefficients $\bnd$, we establish (\ref{eqdw}) by integrating the
above equality over $\theta\in[0,\pi]$ and consolidating
coefficients.\qed
\end{itemize}
\noqed\end{pf}

If $k \geq1$ is an integer, so that $d = 2k + 1$ is odd, the
recursion (\ref{eqdw}) allows us to express the $d$-Schoenberg
coefficient $\bnd$ in terms of the Fourier cosine coefficients
$\bnone, b_{n+2,1},\ldots, b_{n+2k,1}$. For instance,
%
%
\begin{equation}
\label{eqbn5} b_{n,5} = \frac{(n+2) (n+3)}{3 \cdot4} \biggl( \bnone^* - 2
\frac{n+2}{n+3} b_{n+2,1} + \frac{n+1}{n+3} b_{n+4,1} \biggr)
\end{equation}
for all integers $n \geq0$, where $b_{0,1}^* = 2 b_{0,1}$ and
$\bnone^* = \bnone$ for $n \geq1$. Similarly, if $d = 2k + 2$ is
even, we can express the $d$-Schoenberg coefficient $\bnd$ in terms of
the Legendre coefficients $\bntwo, b_{n+2,2},\ldots, b_{n+2k,2}$.
Furthermore, it is possible to relate the Legendre coefficients to the
Fourier cosine coefficients, in that, subject to weak regularity
conditions,
%
%
\begin{equation}
\label{eqbn2} \bntwo= \sum_{k=0}^\infty
c_k^n \bigl( b_{n+2k,1}^* - b_{n+2k+2,1} \bigr),
\end{equation}
where
\[
c_k^n = \frac{2^{2n} (n!)^2}{(2n)!} \frac{1 \cdot3 \cdots(2k-1)}{k!}
\frac{(n+1) \cdot(n+2) \cdots(n+k)}{(2n+3) \cdot(2n+5) \cdots(2n+2k+1)}
\]
for integers $n \geq0$ and $k \geq0$, as used by \citet{Huang2011}.

The recursions in Corollary~\ref{corbnd} allow for dimension walks,
where monotonicity properties of the $d$-Schoenberg coefficients
guarantee positive definiteness in higher dimensions, as follows.

%
\begin{corollary} \label{cordw}
Suppose that the function $\psi: [0,\pi] \to\real$ is continuous
with $\psi(0) = 1$. For integers $d \geq1$ and $n \geq0$,
let $\bnd$ denote the Fourier cosine and Gegenbauer coefficients
(\ref{eqb1}) and (\ref{eqbnd}) of~$\psi$, respectively.
\begin{itemize}[(a)]
\item[(a)]
The function $\psi$ belongs to the class $\Psi_3$ if, and only
if, $b_{2,1} \leq2 b_{0,1}$ and $b_{n+2,1} \leq\bnone$ for
all integers $n \geq1$. It belongs to the class $\Psi_3^+$ if,
and only if, furthermore, the inequality is strict for infinitely many
even and infinitely many odd integers.
\item[(b)]
If $d \geq2$, the function $\psi$ belongs to the class $\Psi_{d+2}$
if, and only if,
%
%
\begin{equation}
\label{eqdwmonotone} b_{n+2,d} \leq\frac{2n+d+3}{2n+d-1}
\frac{(n+d-1) (n+d)}{(n+1) (n+2)} \bnd
\end{equation}
for all integers $n \geq0$. It belongs to the class $\Psi_{d+2}^+$
if, and only if, furthermore, the inequality is strict for infinitely
many even and infinitely many odd integers.
\item[(c)]
If $d \geq2$, the function $\psi$ belongs to the class $\Psi_{d+2}^+$ if
$b_{n+2,d} \leq\bnd$ for all integers $n \geq0$.
\end{itemize}
\end{corollary}


\section{Criteria for positive definiteness and applications} \label
{secapplications}

In this section, we provide easily applicable tests for membership in
the classes $\Psi_d$ and $\Psi_d^+$ and apply them to construct rich
parametric families of correlation functions and radial basis
functions on spheres.

\subsection{Locally supported strictly positive definite functions on
two- and three-dimensional spheres}
\label{seclocally}

There are huge computational savings in the interpolation of scattered
data on spheres if the strictly positive definite function used as the
correlation or radial basis function admits a simple closed form, and
is supported on a spherical cap. Thus, various authors have sought
such functions (\citecs{Schreiner1997}, \citecs{Fasshauer1998},
\citecs{NarcowichWard2002}), with particular emphasis on the
practically most
relevant case of the two-dimensional sphere.

The following result of \citet{Levy1961} constructs a positive
definite function on the circle from a positive definite function on
the real line, in the form of a compactly supported member of the
class $\Phi_1$. Here and in the following, we write
$\varphi_{[0,\pi]}$ for the restriction of a function $\varphi:
[0,\infty) \to\real$ to the interval $[0,\pi]$. For essentially
identical results, see Exercise 1.10.25 of \citet{Sasvari1994} and
\citet{Wood1995}.

%
\begin{theorem} \label{thcompact1}
Suppose that the function $\varphi\in\Phi_1$ is such that
$\varphi(t) = 0 $ for $t \geq\pi$. Then the restriction $\psi=
\varphi_{[0,\pi]}$ belongs to the class $\Psi_1$.
\end{theorem}

On the two-dimensional sphere, a particularly simple way of
constructing a locally supported member of the corresponding class
$\Psi_2^+$ is to take a compactly supported member of the class
$\Phi_3$, such as any of the functions in Table~\ref{tabcompact}, and
apply Yadrenko's recipe (\ref{eqYadrenko}). If $\varphi\in\Phi_3$
is such that $\varphi(t) = 0$ for $t \geq a$, where $a \in(0,\pi)$,
then the function $\psi$ in (\ref{eqYadrenko}) belongs to the class
$\Psi_2^+$ and satisfies $\psi(\theta) = 0$ for $\theta\geq2
\arcsin
\frac{a}{2}$. However, the approach is subject to the problems
described in Section~\ref{secEuclidean}, in that the use of the
chordal distance may result in physically unrealistic representations.

Perhaps surprisingly, the following result shows that compactly
supported, isotropic positive definite functions on the Euclidean
space $\real^3$ remain positive definite with the great circle
distance on a one-, two- or three-dimensional sphere as argument.
Thus, they can be used either with the chordal distance or with the
great circle distance as argument.

%
\begin{theorem} \label{thcompact3}
Suppose that the function $\varphi\in\Phi_3$ is such that
$\varphi(t) = 0$ for $t \geq\pi$. Then the restriction $\psi=
\varphi_{[0,\pi]}$ belongs to the class $\Psi_3^+$.
\end{theorem}

\begin{pf}
If the function $\varphi\in\Phi_3$ is such that $\varphi(t) = 0$
for $t \geq\pi$, the Fourier cosine coefficient (\ref{eqb1}) of
its restriction $\psi= \varphi_{[0,\pi]}$ can be written as
\[
\bnone= \frac{2}{\pi} \int_0^\pi\cos(n
\theta) \psi(\theta) \,\dd\theta= \frac{1}{\pi} \int_\real
\cos(nt) \phi(t) \,\dd t = 2 \hat\phi(n)
\]
for all integers $n \geq1$. Here, the function $\phi: \real\to
\real$ is defined by $\phi(t) = \varphi(|t|)$ for $t \in\real$, and
$\hat\phi(u)= \frac{1} {2 \pi} \int_\real\cos(ut) \phi(t)
\,\dd
t$ denotes the inverse Fourier transform of $\phi$ evaluated at $u \in
\real$. Similarly, $b_{0,1} = \hat\phi(0)$. By equation (36) of
\citet{Gneiting1998b} in concert with the Paley--Wiener theorem
\citep{Paley1934}, the function $\hat\phi(u)$ is strictly decreasing
in $u \geq0$. Thus, the sequence $\bnone$ is strictly decreasing in
$n \geq1$; furthermore, $b_{2,1} = 2 \hat\phi(2) < 2
\hat\phi(0) = 2 b_{0,1}$. By part (a) of Corollary~\ref{cordw},
we conclude that $\psi\in\Psi_3^+$.
\end{pf}

Theorem~\ref{thcompact3} permits us to use any of the functions in
Table~\ref{tabcompact} with support parameter $c \leq\pi$ as a
correlation function or radial basis function with the spherical or
great circle distance on the sphere $\sphere^d$ as argument, where $d
\leq3$. This gives rise to flexible parametric families of locally
supported members of the class $\Psi_d^+$, where $d \leq3$, that
admit particularly simple closed form expressions. A first example
is the \textit{spherical} family,
%
%
\begin{equation}
\label{eqspherical} \psi(\theta) = \biggl( 1 + \frac{1}{2}
\frac{\theta}{c} \biggr) \biggl( 1 - \frac{\theta}{c} \biggr)_+^2
\quad\mbox{for} \quad\theta\in[0,\pi],
\end{equation}
which derives from a popular correlation model in geostatistics.

Using families of compactly supported members of the class $\Phi_3$
described by \citet{Wendland1995} and \citet{Gneiting1999a}, Theorem
\ref{thcompact3} confirms that Askey's (\citeyear{Askey1973})
truncated power function
%
%
\begin{equation}
\label{eqAskey} \psi(\theta) = \biggl( 1 - \frac{\theta}{c}
\biggr)^\tau_+ \quad\mbox{for} \quad\theta\in[0,\pi]
\end{equation}
belongs to the class $\Psi_3^+$ if $\tau\geq2$ and $c \in(0,\pi]$,
as shown by \citet{Beatson2011}. If smoother functions are
desired, the $C^2$-Wendland function
%
%
\begin{equation}
\label{eqWendland15} \psi(\theta) = \biggl( 1 + \tau\frac{\theta}{c}
\biggr) \biggl( 1 - \frac{\theta}{c} \biggr)^\tau_+ \quad\mbox{for}
\quad\theta\in[0,\pi]
\end{equation}
s in $\Psi_3^+$ if $\tau\geq4$ and $c \in(0,\pi]$; similarly, the
$C^4$-Wendland function defined by
%
%
\begin{equation}
\label{eqWendland25} \psi(t) = \biggl( 1 + \tau\frac{\theta}{c} +
\frac{\tau^2-1}{3} \frac{\theta^2}{c^2} \biggr) \biggl( 1 - \frac{\theta}{c}
\biggr)^\tau_+ \quad\mbox{for} \quad\theta\in[0,\pi]
\end{equation}
belongs to the class $\Psi_3^+$ if $\tau\geq6$ and $c \in(0,\pi]$.

\begin{figure}

\includegraphics{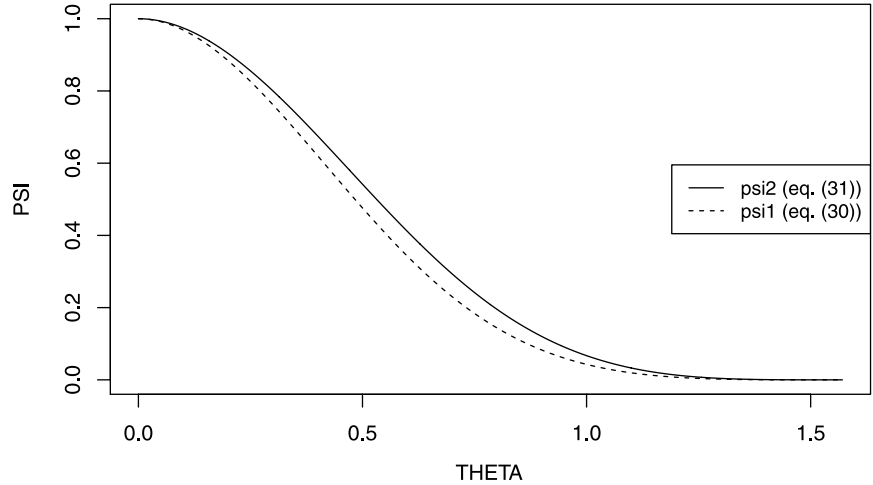}

\caption{The \citet{GaspariCohn1999} localization function with
Euclidean or chordal distance ($\psi_1$; eq.~(\protect\ref
{eqGC1})), and
with spherical or great circle distance ($\psi_2$;
eq. (\protect\ref{eqGC2})), as argument.} \label{figGC}\vspace*{-3pt}
\end{figure}

In atmospheric data assimilation, locally supported positive definite
functions are used for the distance-dependent filtering of spatial
covariance estimates on planet Earth \citep{Hamill2001}. The
traditional construction relies on the \citet{GaspariCohn1999}
function,
\[
\varphi_{\mathrm{GC}}(t) = \cases{ %
1 -
\frac{20}{3} t^2 + 5 t^3 + 8
t^4 - 8 t^5, & $0 \leq t \leq\frac{1}{2}$,
\vspace*{2pt}\cr
\frac{1}{3} t^{-1} \bigl( 8 t^2 + 8 t - 1
\bigr) ( 1 - t )^4, & $\frac{1}{2} \leq t
\leq1$,
\vspace*{2pt}\cr
0, & $t \geq1$,}
\]
which belongs to the class $\Phi_3$, along with all functions of the
form $\varphi(t) = \varphi_{\mathrm{GC}}(t/c)$, where $c > 0$ is a
constant. If $c \in(0,\pi]$, Yadrenko's construction (\ref{eqYadrenko})
yields the localization function
%
%
\begin{equation}
\label{eqGC1} \psi_1(\theta) = \varphi_{\mathrm{GC}} \biggl(
\frac{\sin\frac
{\theta
}{2}}{\sin\frac{c}{2}} \biggr) \quad\mbox{for} \quad\theta\in[0,\pi],
\end{equation}
which is a member of the class $\Psi_2^+$ with support $[0,c]$. Theorem
\ref{thcompact3} suggests a natural alternative, in that, for every
$c \in(0,\pi]$, the function defined by
%
%
\begin{equation}
\label{eqGC2} \psi_2(\theta) = \varphi_{\mathrm{GC}} \biggl(
\frac{\theta}{c} \biggr) \quad\mbox{for} \quad\theta\in[0,\pi]
\end{equation}
is also a member of the class $\Psi_2^+$ with support $[0,c]$.
Clearly, $\psi_2(\theta) > \psi_1(\theta)$ for $\theta\in(0,c)$, as
illustrated in Figure~\ref{figGC}, where $c = \frac{\pi}{2}$. This
suggests that $\psi_2$ might be a more effective localization function
%
%
than the traditional choice, $\psi_1$. Similar comments apply to
covariance tapers in spatial statistics, as proposed by
\citet{Furrer2006}.

\subsection{Criteria of P\'olya type} \label{secPolya}

Criteria of P\'olya type provide simple sufficient conditions for
positive definiteness, by imposing convexity conditions on a candidate
function and/or its higher derivatives. \citet{Gneiting2001a} reviews
and develops P\'olya criteria on Euclidean spaces. In what follows,
we are concerned with analogues on spheres that complement the recent
work of \citet{Beatson2011}.

The following theorem summarizes results of \citet{Wood1995},
\citet{Gneiting1998a} and \citet{Beatson2011} in the case of the
circle $\sphere^1$.

%
\begin{theorem} \label{thPolya1}
Suppose that the function $\psi: [0,\pi] \to\real$ is continuous,
nonincreasing and convex with $\psi(0) = 1$ and $\int_0^\pi
\psi(\theta) \,\dd\theta\geq0$. Then $\psi$ belongs to the
class $\Psi_1$. It belongs to the class $\Psi_1^+$ if it is
not piecewise linear.
\end{theorem}

Our next result is a P\'olya criterion that applies on spheres
$\sphere^d$ of dimension $d \leq3$. Its proof relies on the
subsequent lemmas, the first of which is well known.

%
\begin{theorem} \label{thPolya3}
Suppose that $\varphi: [0,\infty) \to\real$ is a continuous
function with $\varphi(0) = 1$, $\lim_{t \to\infty} \varphi(t) =
0$ and a continuous derivative $\varphi'$ such that $-
\varphi'(t^{1/2})$ is convex for $t > 0$. Then the restriction
$\psi= \varphi_{[0,\pi]}$ belongs to the class $\Psi_3^+$.\vadjust{\goodbreak}
\end{theorem}

To give an example, we may conclude from Theorem~\ref{thPolya3} that
the truncated power function (\ref{eqAskey}) with shape parameter
$\tau\geq2$ belongs to the class $\Psi_3^+$ for all values of the
support parameter $c > 0$, including the case $c > \pi$, in which it
is supported globally.

%
\begin{lemma} \label{lemixture}
Let $d \geq1$ be an integer. Suppose that $C$ is a subset of a
Euclidean space, and let $P$ be a Borel probability measure on
$C$. If $\psi_c$ belongs to the class $\Psi_d $ for every $c
\in C$, then the function $\psi$ defined by
\[
\psi(\theta) = \int_C \psi_c(\theta) \,\dd
P(c) \quad\mbox{for} \quad\theta\in[0,\pi]
\]
belongs to the class $\Psi_d$, too. If furthermore $\psi_c$
belongs to $\Psi_d^+$ for every $c$ in a set of positive
$P$-measure, then $\psi$ belongs to the class $\Psi_d^+$.
\end{lemma}

%
\begin{lemma} \label{lespherical}
For all $c > 0$, the function defined by
%
%
\begin{equation}
\label{eqsphericalpsi} \psi_c(\theta) = \biggl( 1 +
\frac{1}{2} \frac{\theta}{c} \biggr) \biggl( 1 - \frac{\theta}{c}
\biggr)_+^2 \quad\mbox{for} \quad\theta\in[0,\pi]
\end{equation}
belongs to the class $\Psi_3^+$.
\end{lemma}

\begin{pf}
By Theorem~\ref{thcompact3}, the function $\psi_c$ belongs to the
class $\Psi_3^+$ if $c \leq\pi$. If $c \geq\pi$, the formulas in
(\ref{eqb1}) yield the Fourier cosine coefficients
\[
b_{0,1} = \frac{1}{8c^3} \bigl(8c^3 - 6\pi
c^2 + \pi^3\bigr) \qquad\mbox{and} \qquad
b_{2k,1} = \frac{3\pi}{4c^3} \frac{1}{k^2}
\]
for all integers $k \geq1$. Similarly, we find that
\[
b_{2k+1,1} = \frac{3}{\pi c^3} \biggl( \frac{2c^2 - \pi^2}{(2k+1)^2} +
\frac{4}{(2k+1)^4} \biggr)
\]
for all integers $k \geq0$ and conclude from part (a) of Corollary
\ref{cordw} that $\psi_c$ belongs to $\Psi_3^+$.
\end{pf}

\begin{pf*}{Proof of Theorem~\ref{thPolya3}}
By Theorem 3.1 of
\citet{Gneiting1999b}, the function $\varphi$ admits a representation
of the form
\[
\varphi(t) = \int_{(0,\infty)} \varphi_c(t) \,\dd F(c)
\quad\mbox{for} \quad t \geq0,
\]
where $\varphi_c$ is defined by
%
%
\begin{equation}
\label{eqsphericalEuclidean} \varphi_c(t) = \biggl( 1 +
\frac{1}{2} \frac{t}{c} \biggr) \biggl( 1 - \frac{t}{c}
\biggr)_+^2 \quad\mbox{for} \quad t \geq0,
\end{equation}
and $F$ is a probability measure on $(0,\infty)$. Therefore, the
restriction $\psi= \varphi_{[0,\pi]}$ is of the form
\[
\psi(\theta) = \int_{(0,\infty)} \psi_c(\theta) \,\dd
F(c) \quad\mbox{for} \quad\theta\in[0,\pi],
\]
where $\psi_c$ is defined by (\ref{eqsphericalpsi}). By Lemma
\ref{lespherical}, the function $\psi_c$ belongs to the class
$\Psi_3^+$ for all $c > 0$, and so we conclude from Lemma
\ref{lemixture} that $\psi$ is in $\Psi_3^+$.
\end{pf*}

Our next theorem is a slight generalization of the key result in
Beatson, zu~Castell and
Xu (\citeyear{Beatson2011}). In analogy to the corresponding results on
Euclidean spaces \citep{Gneiting1999b}, the case $n = 1$ yields a more
concise, but less general, criterion than Theorem~\ref{thPolya3}.

%
\begin{theorem} \label{thPolyan}
Let $n \leq3$ be a positive integer. Suppose that $\varphi:
[0,\infty) \to\real$ is a continuous function with $\varphi(0) =
1$, $\lim_{t \to\infty} \varphi(t) = 0$ and a derivative
$\varphi^{(n)}$ of order $n$ such that $(-1)^n \varphi^{(n)}(t)$
is convex for $t > 0$. Then the restriction $\psi=
\varphi_{[0,\pi]}$ belongs to $\Psi_{2n+1}^+$.
\end{theorem}

In a recent \textit{tour de force}, \citet{Beatson2011} demonstrate the
following remarkable result, which is the key to proving Theorem
\ref{thPolyan}.

%
\begin{lemma} \label{leBCX}
Let $n \leq3$ be a positive integer, and suppose that $c \in(0,\pi]$.
Then the function defined by
%
%
\begin{equation}
\label{eqtpower} \psi_{c,n}(\theta) = \biggl( 1 - \frac{\theta}{c}
\biggr)_+^{n+1} \quad\mbox{for} \quad\theta\in[0,\pi]
\end{equation}
belongs to the class $\Psi_{2n+1}^+$.
\end{lemma}

The next lemma concerns the case $c > \pi$, in which the truncated
power function (\ref{eqtpower}) is supported globally.

%
\begin{lemma} \label{letpower}
For all integers $n \geq1$ and all real numbers $c \geq\pi$, the
function $\psi_{c,n}$ in (\ref{eqtpower}) belongs to the
class $\Psi_\infty^+$.
\end{lemma}

\begin{pf}
By a convolution argument due to Fuglede (\citecs{Hjorth1998}, p. 272),
the function $\psi_{\pi,0}(\theta) = 1 - \frac{\theta}{\pi}$ belongs
to the class $\Psi_\infty$. As the class $\Psi_\infty$ is convex and
closed under products, we see that if $n \geq0$ is an integer and $c
\geq\pi$ then
\[
\psi_{c,n}(\theta) = \biggl( \biggl( 1 - \frac{\pi}{c} \biggr) +
\frac{\pi}{c} \psi_{\pi,0}(\theta) \biggr)^{n+1} \quad
\mbox{for} \quad\theta\in[0,\pi],
\]
whence $\psi_{c,n}$ is in the class $\Psi_\infty$, too. By a
straightforward direct calculation, the Fourier cosine coefficients of
$\psi_{c,1}$ and $\psi_{c,2}$ are strictly positive for all $c \geq
\pi$. Therefore, by Theorem 1 of \citet{Xu1992} and the fact that the
class $\Psi_1^+$ is closed under products, the function $\psi_{c,n}$
is in the class $\Psi_1^+$ for all integers $n \geq1$ and all $c \geq
\pi$. Using part (a) of Corollary~\ref{corPsi}, we conclude that
$\psi_{c,n} \in\Psi_\infty^+$ for all integers $n \geq1$ and all $c
\geq\pi$.
\end{pf}

\begin{pf*}{Proof of Theorem~\ref{thPolyan}} By Theorem 3.1 of
\citet{Gneiting1999b}, the function $\varphi$ admits a representation
of the form
\[
\varphi(t) = \int_{(0,\infty)} \varphi_{c,n}(t) \,\dd F(c)
\quad\mbox{for} \quad t \geq0,
\]
where the function $\varphi_{c,n}: [0,\infty) \to\real$ is defined
by $\varphi_{c,n}(t) = ( 1 - \frac{\theta}{c} )^{n+1}$
for $t \geq0$ and $F$ is a probability measure on $(0,\infty)$.
Therefore, the restriction $\psi= \varphi_{[0,\pi]}$ is of the form
\[
\psi(\theta) = \int_{(0,\infty)} \psi_{c,n}(\theta) \,\dd
F(c) \quad\mbox{for} \quad\theta\in[0,\pi],
\]
where $\psi_{c,n}$ is given by (\ref{eqtpower}). By Lemmas
\ref{leBCX} and~\ref{letpower}, the function $\psi_{c,n}$ belongs to
the class $\Psi_{2n+1}^+$ for all $c > 0$. Therefore, we
conclude from Lemma~\ref{lemixture} that $\psi$ is in
$\Psi_{2n+1}^+$, too.
\end{pf*}

\citet{Beatson2011} conjectured that the statement of Lemma
\ref{leBCX} holds for all integers $n \geq1$. In view of Lemma
\ref{letpower}, if their conjecture is true, the statement of Theorem
\ref{thPolyan} holds for all integers $n \geq1$, too, with the
proof being unchanged.

It is interesting to observe that Theorem~\ref{thPolyan} is a
stronger result than Theorem 1.3 in \citet{Beatson2011}, in that it
does not impose any support conditions on the candidate function
$\psi$. In contrast, the latter criterion requires $\psi$ to be
locally supported. Such an assumption also needs to made in the
spherical analogue of Theorem 1.1 in \citet{Gneiting2001a}, which we
state and prove in Appendix A in the supplemental article
\citep{Gneiting2012}.

\subsection{Completely monotone functions} \label{seccm}

A function $\varphi: [0,\infty) \to\real$ is completely monotone if
it possesses derivatives $\varphi^{(k)}$ of all orders with $(-1)^k
\varphi^{(k)}(t) \geq0$ for all integers $k \geq0$ and all $t > 0$.
Our next result shows that the restrictions of
completely monotone functions belong to the class $\Psi_\infty^+$.

%
\begin{theorem} \label{thcm}
Suppose that the function $\varphi: [0,\infty) \to\real$ is
completely monotone with $\varphi(0) = 1$ and not constant.
Then the restriction $\psi= \varphi_{[0,\pi]}$ belongs to the
class $\Psi_\infty^+$.
\end{theorem}

\begin{pf}
Let $a > 0$ and consider the truncated power function
(\ref{eqtpower}). By Lemma~\ref{letpower}, $\psi_{n/a,n}$ belongs
to the class $\Psi_\infty^+$ for all sufficiently large integers
$n$. Hence, the function $\psi_a$ defined by
\[
\psi_a(\theta) = e^{- a \theta} = \lim_{n \to\infty}
\psi_{n/a,n}(\theta) \quad\mbox{for} \quad\theta\in[0,\pi]
\]
belongs to the class $\Psi_\infty$, too. By Theorem~\ref{thPolya3},
$\psi_a \in\Psi_3^+$, which in view of part (a) of Corollary
\ref{corPsi} implies that $\psi_a \in\Psi_\infty^+$. Now suppose
that the function $\varphi: [0,\infty) \to\real$ satisfies the
conditions of the theorem. Invoking Bernstein's Theorem, we see
that $\varphi$ admits a representation of the form
\[
\varphi(t) = \int_{[0,\infty)} e^{-at} \,\dd F(a) \quad
\mbox{for} \quad t \geq0,
\]
where $F$ is a probability measure on $[0,\infty)$. Therefore, the
restriction $\psi= \varphi_{[0,\pi]}$ is of the form
\[
\psi(\theta) = \int_{[0,\infty)} e^{-a \theta} \,\dd F(a) = \int
_{[0,\infty)} \psi_a(\theta) \,\dd F(a) \quad\mbox{for}
\quad\theta\in[0,\pi].
\]
Since $\varphi$ is not constant, $F$ has mass away from the origin.
As $\psi_a$ is in the class $\Psi_\infty^+$ for all $a > 0$, we
conclude from Lemma~\ref{lemixture} that $\psi$ belongs to the class
$\Psi_\infty^+$.
\end{pf}

\citet{Miller2001} present a wealth of examples of completely monotone
functions, which thus can serve as correlation functions or radial
basis functions on spheres of any dimension.

\subsection{Necessary conditions} \label{secnecessary}

The next result states conditions that prohibit membership in the
class $\Psi_1$, and therefore in any of the classes $\Psi_d$.
Heuristically, the common motif can be paraphrased as follows:

\begin{quote}
If a positive definite function admits a certain degree of smoothness
at the origin, it admits the same degree of smoothness everywhere.
\end{quote}

The conditions in parts (a) through (d) of the following theorem can
be interpreted as formal descriptions of circumstances under which
this overarching principle is violated. Part (a) corresponds to a
well-known result in the theory of characteristic functions, parts (b)
and (c) are due to \citet{Wood1995} and \citet{Gneiting1998a},
respectively, and part (d) rests on a result in \citet{Devinatz1959}.

%
\begin{theorem} \label{thnecessary}
Suppose that the function $\psi: [0,\pi] \to\real$ can be
represented as the restriction $\psi= \phi_{[0,\pi]}$ of an even
and continuous function $\phi: \real\to\real$ that satisfies any
of the following conditions.
\begin{itemize}[(a)]
\item[(a)] For some integer $k \geq1$ the derivative
$\phi^{(2k)}(0)$ exists, but $\phi$ fails to be $2k$
times differentiable on $(0,\pi)$.
\item[(b)] For some integer $k \geq1$, the function $\phi$ is
$2k - 1$ times continuously differentiable on $[-\pi,\pi]$ with
$\phi^{(2k-1)}(\pi) \not= 0$ and $|\phi^{(2k-1)}(t)| \leq b_1
|t|^{\beta_1}$ for some $b_1 > 0$ and $\beta_1 > 0$, and is $2k +
1$ times continuously differentiable on $[-\pi,0) \cup(0,\pi]$
with $|\phi^{(2k+1)}(t)| \leq b_2 |t|^{-\beta_2}$ for some $b_2 > 0$
and $\beta_2 \in(1,2)$.
\item[(c)] For some integer $k \geq1$, the function $\phi$ is
$2k$ times differentiable on $[-\pi,\pi]$ with $\phi^{(j)}(\pi)
\not= 0$ for some odd integer $j$, where $1 \leq j \leq2k-1$.
\item[(d)] The function $\phi$ is analytic and not of period
$2\pi$.
\end{itemize}
Then $\psi$ does not belong to the class $\Psi_1$.
\end{theorem}

An important caveat is the possibly surprising fact that, given some
function $\varphi: [0,\infty) \to\real$, the mapping $\psi_c:
[0,\pi] \to\real$ defined by
\[
\psi_c(\theta) = \varphi\biggl( \frac{\theta}{c} \biggr) \quad
\mbox{for} \quad\theta\in[0,\pi]
\]
may belong to the class $\Psi_d$ for some specific values of the scale
parameter $c > 0$ only, rather than for all $c > 0$. For example, the
$C^2$- and $C^4$-Wendland functions in Table~\ref{tabPsi2} allow for
$c \in(0,\pi]$ only, with part (c) of Theorem~\ref{thnecessary}
excluding the case $c > \pi$.

A notable exception occurs when $\varphi$ is completely monotone with
$\varphi(0) = 1$ and not constant. Then the mapping $t \mapsto
\varphi(t/c)$ is completely monotone for all $c > 0$, and we conclude
from Theorem~\ref{thcm} that $\psi_c \in\Psi_\infty^+$ for all
values of the scale parameter $c > 0$. The first four entries in
Table~\ref{tabPsi2} are of this type and will be discussed below.

\subsection{Examples} \label{secexamples}

We now apply the criteria of Sections~\ref{seccm} and
\ref{secnecessary} to study parametric families of globally supported
correlation functions and radial basis functions on spheres. In doing
so, we think of $c > 0$ as a scale or support parameter, $\alpha> 0$
or $\nu> 0$ as a smoothness parameter and $\tau> 0$ as a shape
parameter.

We consider the families in Table~\ref{tabPhiinfty}, namely the
powered exponential, Mat\'ern, generalized Cauchy and Dagum classes.
For the parameter values stated in Table~\ref{tabPhiinfty} these
families serve as isotropic correlation functions on Euclidean spaces,
where the argument is the Euclidean distance. Investigating whether
they can serve as isotropic correlation functions with spherical or
great circle distance on $\sphere^d$ as argument, we supplement and
complete results of \citet{Huang2011} in the case $d = 2$, as
summarized in the first four entries of Table~\ref{tabPsi2}. In
Examples~\ref{exexponential} and~\ref{exMatern}, the results and
proofs do not depend on the scale parameter $c > 0$, and so we only
discuss the smoothness parameter $\alpha$ or $\nu$, respectively.

%
\begin{example}[(Powered exponential family)] \label{exexponential}
The members of the powered exponential family are of the form
%
%
\begin{equation}
\label{eqexponential} \psi(\theta) = \exp\biggl( - \biggl( \frac{\theta}{c}
\biggr)^{
\alpha} \biggr) \quad\mbox{for} \quad\theta\in[0,\pi].
\end{equation}
If $\alpha\in(0,1]$, a straightforward application of Theorem
\ref{thcm} shows that $\psi\in\Psi_\infty^+$. If $\alpha> 1$, we
see from parts (b) and (c) of Theorem~\ref{thnecessary} that $\psi$
does not belong to the class $\Psi_1$, where part (b) applies when
$\alpha\in(1,2)$, and part (c) when $\alpha\geq2$, both using $k =
1$. Therefore, if $\alpha> 1$ then $\psi$ does not belong to any of
the classes $\Psi_d$.
\end{example}

%
\begin{example}[(Mat\'ern family)] \label{exMatern}
The members of the Mat\'ern family can be written as
%
%
\begin{equation}
\label{eqMatern} \psi(\theta) = \frac{2^{\nu-1}}{\Gamma(\nu)} \biggl(
\frac{\theta}{c}
\biggr)^{ \nu} K_\nu\biggl( \frac
{\theta
}{c} \biggr) \quad
\mbox{for} \quad\theta\in[0,\pi],
\end{equation}
where $K_\nu$ denotes the modified Bessel function of the second kind
of order $\nu$. If $\nu= n + \frac{1}{2}$, where $n \geq0$ is an
integer, then $\psi$ equals the product of $\exp(-\theta/c)$ and a
polynomial of degree $n$ in $\theta$
\citep{GuttorpGneiting2006}.

If $\nu\in(0,\frac{1}{2}]$, then $\psi$ belongs to the class
$\Psi_\infty^+$. Indeed, by Theorem 5 of \citet{Miller2001} the
function $\phi_1(t) = \exp(t) t^\nu K_\nu(t)$ is completely
monotone if $\nu\leq\frac{1}{2}$. By Theorem 1 of
\citet{Miller2001} the function $\phi_0(t) = \exp(-t) \phi_1(t) =
t^\nu K_\nu(t)$ is completely monotone, too. In view of Theorem
\ref{thcm}, this proves our claim.

If $\nu> \frac{1}{2}$, then $\psi$ does not belong to the class
$\Psi_1$, and thus neither to any of the classes $\Psi_d$. If $\nu
\in(\frac{1}{2},1)$, this follows from an application of the
estimates in Section 10.31 of \citet{DLMF} to part (b) of Theorem
\ref{thnecessary}, where $k = 1$, $\beta_1 = 2\nu- 1$ and $\beta_2 =
2\nu- 3$. If $\nu\geq1$, we apply part (c) of Theorem
\ref{thnecessary}, where $k = 1$.
\end{example}

In our next example, the results do not depend on the scale parameter
$c > 0$ and the shape parameter $\tau> 0$, and so we discuss the
smoothness parameter $\alpha$ only.

%
\begin{example}[(Generalized Cauchy family)] \label{exCauchy}
The members of the generalized Cauchy family are of the form
%
%
\begin{equation}
\label{eqCauchy} \psi(\theta) = \biggl( 1 + \biggl( \frac{\theta}{c}
\biggr)^{ \alpha} \biggr)^{-
\tau
/\alpha} \quad\mbox{for} \quad\theta\in[0,
\pi].
\end{equation}
If $\alpha\in(0,1]$, a straightforward application of Theorem
\ref{thcm} shows that $\psi$ is a member of the class~$\Psi_\infty^+$.
However, if $\alpha> 1$ the function $\psi$ in (\ref{eqCauchy}) does
not belong to the class $\Psi_1$. If $\alpha \in(1,2)$, this is evident
from part (b) of Theorem~\ref{thnecessary}, where $k = 1$, $\beta_1 =
\alpha- 1$ and $\beta_2 = \alpha- 3$. If $\alpha\geq2$ we apply part
(c) of Theorem~\ref{thnecessary}, where $k = 1$.
\end{example}

%
\begin{example}[(Dagum family)] \label{exDagum}
The members of the Dagum family are of the form
%
%
\begin{equation}
\label{eqDagum} \psi(\theta) = 1 - \biggl( \biggl( \frac{\theta}{c}
\biggr)^{ \tau} \bigg/ \biggl( 1 + \frac
{\theta}{c} \biggr)^{ \tau}
\biggr)^{\alpha/\tau} \quad\mbox{for} \quad\theta\in[0,\pi].
\end{equation}
Using a comment by Berg, Mateu and Porcu (\citeyear{Berg2008}, p. 1146)
along with Theorem~\ref{thcm}, we see that if $c > 0$, $\tau\in(0,1]$
and $\alpha\in(0,\tau)$ then $\psi$ is a member of the class
$\Psi_\infty^+$.
\end{example}

Another example that involves oscillating trigonometric and Bessel
functions is given in Appendix B in the supplemental article
\citep{Gneiting2012}.

In the context of isotropic random fields and random particles, the
smoothness properties of the associated random surface are governed by
the behavior of the correlation function at the origin. Specifically,
if a function $\psi\in\Psi_2$ admits the relationship
(\ref{eqfractalpsi}) for some $\alpha\in(0,2]$, the corresponding
Gaussian sample paths on the two-dimensional sphere have fractal or
Hausdorff dimension $3 - \frac{\alpha}{2}$ almost surely
\citep{Hansen2011}. In this sense, the results in the above examples
are restrictive, in that the smoothness parameter needs to satisfy
$\alpha\leq1$ or $\nu\leq\frac{1}{2}$, respectively. Parametric
families of correlation\vadjust{\goodbreak} functions that admit the full range of viable
Hausdorff dimensions include the sine power family
(\ref{eqsinepower}) and the convolution construction in Section 4.3
of \citet{Hansen2011}. Alternatively, Yadrenko's construction
(\ref{eqYadrenko}) can be applied to any of the parametric families
in Table~\ref{tabPhiinfty}.


\section{Challenges for future work} \label{secdiscussion}

In this paper, we have reviewed and developed characterizations and
constructions of isotropic positive definite functions on spheres, and
we have applied them to provide rich parametric classes of such
functions. Some of the key results in the practically most relevant
cases of one-, two- and three-dimensional spheres are summarized in
Table~\ref{tabPsi2}. Whenever required, the closure properties of
the classes of the positive definite or strictly positive definite
functions offer additional flexibility. For example, while all
entries in Table~\ref{tabPsi2} yield nonnegative correlations only,
we can easily model negative correlations, by using convex sums or
products of an entry in the table with a suitable Gegenbauer function,
such as $\psi(\theta) = \cos\theta$, which belongs to the class
$\Psi_\infty^-$.

Despite substantial advances in the study of positive definite
functions on spheres, many interesting and important questions remain
open. Therefore, Appendix C in the supplemental article
\citep{Gneiting2012} describes 18 research problems that aim to
stimulate future research in mathematical analysis, probability theory
and spatial statistics. The problems vary in scope and difficulty,
range from harmonic analysis to statistical methodology, and include
what appear to be tedious but routine questions, such as Problem 1,
along with well known major challenges, such as Problems 14 through
16, which have been under scrutiny for decades. The first eight
problems are of an analytic character and concern the characterization
and breadth of the classes $\Psi_d$ and~$\Psi_d^+$. Then we state
open questions about the parameter spaces for various types of
correlation models, relate to the fractal index and the sample path
properties of Gaussian random fields, and turn to challenges in
spatial and spatio-temporal statistics.

\section*{Acknowledgements}

The author thanks Elena Berdysheva, Wolfgang zu Castell, Werner Ehm,
Peter Guttorp, Linda Hansen, Matt\-hias Katzfuss, Finn Lindgren,
Emilio Porcu, Zolt\'an Sasv\'ari, Michael Scheuerer, Thordis
Thorarinsdottir, Jon Wellner, Johanna Ziegel, two anonymous referees
and the editor, Richard Davis, for comments and discussions. His
research has been supported by the Alfried Krupp von Bohlen und
Halbach Foundation, and by the German Research Foundation within the
programme ``Spatio-/Temporal Graphical Models and Applications in
Image Analysis'', grant GRK 1653.

\begin{supplement}
\stitle{Supplement to ``Strictly and non-strictly positive definite
functions on spheres''}
\slink[doi]{10.3150/12-BEJSP06SUPP} 
\sdatatype{.pdf}
\sfilename{BEJSP06\_supp.pdf}
\sdescription{Appendix A states and proves further criteria of P\'olya
type, thereby complementing Section~\ref{secPolya}. Appendix B
studies an example that involves oscillating trigonometric and
Bessel functions, as hinted at in Section~\ref{secexamples}.
Appendix C describes open problems that aim to stimulate future
research in areas ranging from harmonic analysis to spatial
statistics.}
\end{supplement}

%

\printhistory

\end{document}